\documentclass[11pt,reqno]{amsart}
\usepackage[T1]{fontenc}
\usepackage{lmodern}
\usepackage{amsmath,amssymb}
\usepackage[expansion=false]{microtype}
\usepackage{needspace}
\usepackage{tikz}
\usepackage[a4paper,hmargin=28mm,vmargin=25mm]{geometry}
\usepackage{xurl}
\usepackage{hyperref}

\newtheorem{theorem}{Theorem}[section]
\newtheorem{proposition}[theorem]{Proposition}
\newtheorem{lemma}[theorem]{Lemma}
\newtheorem{corollary}[theorem]{Corollary}
\theoremstyle{remark}
\newtheorem{example}[theorem]{Example}
\newtheorem{remark}[theorem]{Remark}
\numberwithin{equation}{section}
\title[Tilings of prime-power integer tiles]{Tilings of integer sets\\
of prime-power cardinality}
\author[W. Hu]{Wenchuan Hu}
\address[Wenchuan Hu]{Sichuan University}
\email{wenchuan@scu.edu.cn}
\author[C. Wang]{Chao Wang}
\address[Chao Wang]{Nankai University}
\email{wangchao@nankai.edu.cn}
\author[E. Wang]{Erxiao Wang}
\address[Erxiao Wang]{Zhejiang Normal University}
\email{wang.eric@zjnu.edu.cn}
\author[F. Zhu]{Fuhai Zhu}
\address[Fuhai Zhu]{Nanjing University}
\email{zhufuhai@nju.edu.cn}
\thanks{Erxiao Wang is the corresponding author. Research was supported by
the National Natural Science Foundation of China (NSFC-RGC 12361161603).}
\date{}

\begin{document}
\begin{abstract}
Let $A$ be a finite integer tile of cardinality $p^k$, where $p$ is prime
and $k\geq1$. We first determine all possible fundamental periods of
tilings of $A$. A cyclotomic argument gives a common period $PD$, where
$P$ is a power of $p$ and $D$ is a greatest common divisor of normalized
differences in $A$. Each $Pd$ with $d\mid D$ occurs as a fundamental period.
At modulus $Pd$, a chain of greatest common divisors controls the
compatibility of the Chinese remainder fibres. This yields independent
digit functions and exact counts of the tilings. For tilings of fundamental
period $Pd$, we classify their partitions into translates modulo $Pd$ and
count the tilings that induce each partition. The same arithmetic data classify
tiles with identical tiling families. Each class has a unique normalized
tile of minimum diameter whose mask polynomial divides the masks of all
normalized tiles in the class.
\end{abstract}
\maketitle

\section{Introduction}

For subsets $U,V$ of an abelian group, $U\oplus V$ denotes their
sumset when every element of it has a unique representation as $u+v$
with $u\in U$ and $v\in V$. A finite nonempty set $A\subset\mathbb{Z}$
is a tile if $A\oplus T=\mathbb{Z}$ for some $T\subset\mathbb{Z}$.
We call such a set $T$ a \emph{tiling} of $A$. We classify all tilings
of a fixed tile of prime-power cardinality.

Throughout, $A$ is a tile with $|A|=p^k$, where $p$ is prime and $k\geq1$.
We translate $A$ so that $\min A=0$, which preserves its tiling family.
Let $v_p(a)$ be the exponent of $p$ in a nonzero integer $a$.
Newman's criterion \cite[Theorem~1]{Newman}, together with the branching
count recalled in Section~\ref{sec:digits}, gives
\begin{equation}\label{eq:levels}
\{v_p(a-a'):a,a'\in A,\ a\neq a'\}
=\{\nu_1-1,\ldots,\nu_k-1\},
\qquad 1\leq\nu_1<\cdots<\nu_k.
\end{equation}
Define
\begin{equation}\label{eq:data}
n=\nu_k,\qquad P=p^n,\qquad Q_i=p^{\nu_i-1},
\end{equation}
and define
\begin{equation}\label{eq:gcd-data}
\Delta_i=\gcd\left\{\frac{|a-a'|}{Q_i}:
a,a'\in A,\ a\neq a',\ a\equiv a'\pmod{Q_i}\right\},
\qquad D=\Delta_k.
\end{equation}
The integers $\Delta_i$ are well defined and satisfy
\begin{equation}\label{eq:chain}
p\nmid\Delta_i,\qquad \Delta_1\mid\cdots\mid\Delta_k.
\end{equation}

Every tiling is periodic \cite[Lemma~1.2]{CovenMeyerowitz}. Its
\emph{fundamental period} is $\operatorname{per}(T)=\min\{N\geq1:T+N=T\}$.
Let $\mathcal{M}(A)$ be the set of these periods over all tilings of $A$.
For a finite set $U$ of nonnegative integers, its mask polynomial is
$U(x)=\sum_{u\in U}x^u$. We write $\Phi_s(x)$ for the $s$th
cyclotomic polynomial and take polynomial divisibility in $\mathbb{Z}[x]$.

Our first theorem gives the possible fundamental periods. This reduces
the classification to the finite cyclic groups of order $Pd$ with $d\mid D$.

\begin{theorem}\label{thm:periods}
The possible fundamental periods of tilings of $A$ are given by
\begin{equation}\label{eq:period-set}
\mathcal{M}(A)=\{Pd:d\mid D\}=\{M\geq1:P\mid M,\ \Phi_M(x)\mid A(x)\}.
\end{equation}
In particular, $PD$ is the largest fundamental period and the least
common period of all tilings of $A$.
\end{theorem}

For $N\geq1$, let $\mathbb{Z}_N=\mathbb{Z}/N\mathbb{Z}$ and let $r_N$
denote reduction modulo $N$, also from $\mathbb{Z}_L$ when $N\mid L$.
For $B\subset\mathbb{Z}_N$, write $B+N\mathbb{Z}$ for its inverse image
in $\mathbb{Z}$. The notation $A\oplus B=\mathbb{Z}_N$ means that
addition modulo $N$ gives a bijection from $A\times B$ onto
$\mathbb{Z}_N$. Equivalently, $B+N\mathbb{Z}$ is a tiling of $A$
\cite[Lemma~1.3]{CovenMeyerowitz}.
We also call such a residue set $B$ a tiling of $A$ in $\mathbb{Z}_N$.

For $d\mid D$, define $d_i=\gcd(d,\Delta_i)$. Then
$d_1\mid\cdots\mid d_k=d$. Let
\[
\rho_d:\mathbb{Z}_{Pd}\longrightarrow\mathbb{Z}_P\times\mathbb{Z}_d,
\qquad \rho_d(b)=(r_P(b),r_d(b))
\]
be the Chinese remainder isomorphism. Write $\overline{A}=r_P(A)$ and let
$\pi_i:\mathbb{Z}_P\longrightarrow\mathbb{Z}_{p^{\nu_i}}$ be reduction.
For $B\subset\mathbb{Z}_{Pd}$, define $C_y\subset\mathbb{Z}_P$ by
\[
\rho_d(B)=\bigcup_{y\in\mathbb{Z}_d}(C_y\times\{y\}).
\]

The classification combines a tiling condition in each fibre with a
compatibility condition on projections of different fibres.

\begin{theorem}\label{thm:classification}
Let $d\mid D$ and $B\subset\mathbb{Z}_{Pd}$. Then
$A\oplus(B+Pd\mathbb{Z})=\mathbb{Z}$ if and only if
\begin{align}
\overline{A}\oplus C_y&=\mathbb{Z}_P
&& (y\in\mathbb{Z}_d),\label{eq:fibre-tiling}\\
y\equiv z\pmod{d_i}\quad&\Longrightarrow\quad
\pi_i(C_y)=\pi_i(C_z)
&& (1\leq i\leq k,\ y,z\in\mathbb{Z}_d).\label{eq:compatibility}
\end{align}
It suffices to impose \eqref{eq:compatibility} for
$i<k$ with $d_i<d_{i+1}$.
For such a tiling, let
\[
H_B=\{h\in\mathbb{Z}_d:C_{y+h}=C_y\text{ for every }y\in\mathbb{Z}_d\}.
\]
Its fundamental period is
\begin{equation}\label{eq:fibre-fundamental-period}
\operatorname{per}(B+Pd\mathbb{Z})=\frac{Pd}{|H_B|}.
\end{equation}
\end{theorem}

Among these tilings, those with fundamental period $Pd$ are exactly
the ones with trivial stabilizer $H_B$.

Two tiles are \emph{tiling-equivalent} if they have the same tilings.
For a finite nonempty set $U\subset\mathbb{Z}$, write
$\operatorname{diam}(U)=\max U-\min U$. The next theorem gives a
canonical representative for each tiling-equivalence class.

\begin{theorem}\label{thm:canonical}
The set
\begin{equation}\label{eq:canonical}
A^*=\bigoplus_{i=1}^k Q_i\Delta_i\{0,\ldots,p-1\}
\end{equation}
is tiling-equivalent to $A$. Every integer tile tiling-equivalent to $A$
has cardinality $p^k$. Two tiles of this cardinality are tiling-equivalent
if and only if their data
\[
\bigl(p,k,(\nu_i)_{i=1}^k,(\Delta_i)_{i=1}^k\bigr)
\]
agree. Moreover,
\begin{equation}\label{eq:canonical-divisibility}
A^*(x)=\prod_{i=1}^k\Phi_p(x^{Q_i\Delta_i})\mid A(x),
\qquad
(p-1)\sum_{i=1}^k Q_i\Delta_i\leq\operatorname{diam}(A).
\end{equation}
Equality in the diameter inequality holds precisely when $A=A^*$.
Consequently, $A^*$ is the unique normalized tile of minimum diameter
in its tiling-equivalence class.
\end{theorem}

For tiles of prime-power cardinality, the construction of
Coven--Meyerowitz gives a common tiling whenever their prime-power
cyclotomic orders agree \cite[proof of Theorem~A and following
remarks]{CovenMeyerowitz}. Those orders determine the smallest possible
fundamental period. The full gcd chain determines the entire tiling
family. Greenfeld--Tao obtain universal periods for
one-dimensional tilings \cite[Corollary~3.5 and Remark~3.6]{GreenfeldTao}.
Our results give the exact least common period in the prime-power case.
The distinction between the fundamental period of a given tiling and
the smallest one over the family is discussed in
\cite[Definition~1]{LabaZakharov}.

Sections~\ref{sec:periods}, \ref{sec:classification}, and
\ref{sec:canonical} prove the three theorems above, respectively.
Section~\ref{sec:counting} counts tilings, Section~\ref{sec:partitions}
classifies and counts their partitions into translates modulo a fundamental
period, and Section~\ref{sec:affine} counts orbits of tilings under affine maps.

\section{Preliminary lemmas}\label{sec:digits}

We recall the standard prime-power structure and dilation results, then
establish the difference restrictions needed for the classification.
Number the base-$p$ digit positions starting with $1$ for the coefficient
of $p^0$, and set $I=\{\nu_1,\ldots,\nu_k\}$ and
$J=\{1,\ldots,n\}\setminus I$. At position $j$, an
occupied class modulo $p^{j-1}$ has \emph{full branching} if all $p$
lifts modulo $p^j$ occur and \emph{single branching} if exactly one occurs.

Newman's criterion bounds the number of distinct valuations of nonzero
differences by $k$. If there are $h$ such valuations, the residue tree
can branch at only $h$ positions, with at most $p$ lifts at each, so
$|A|\leq p^h$. Since $|A|=p^k$, there are exactly $k$ branching
positions, as in \eqref{eq:levels}. Reduction modulo $P$ is injective
because every nonzero difference has valuation less than $n$.
The same cardinality bound forces full branching at every occupied
class at each position in $I$, and single branching elsewhere.
This is the homogeneous-tree structure of \cite[Theorem~4.2]{FanFanShi}.
At level $\nu_i$, the image of $A$ has $p^i$ points, each with
$p^{k-i}$ preimages. In particular, every occupied class modulo $Q_i$
has all $p$ lifts modulo $pQ_i$.
Coven--Meyerowitz describe the recursion in the pure prime-power case
\cite[Section~4]{CovenMeyerowitz}.

We use Tijdeman's dilation theorem \cite{Tijdeman}
in the finite cyclic form of Coven--Meyerowitz
\cite[Lemma~3.1]{CovenMeyerowitz}.

\begin{lemma}\label{lem:dilation}
If $U\oplus V=\mathbb{Z}_N$ and $\lambda\geq1$ satisfies
$\gcd(\lambda,|U|)=1$, then
\begin{equation}\label{eq:dilation}
\lambda U\oplus V=\mathbb{Z}_N.
\end{equation}
The map $u\mapsto\lambda u$ is injective on $U$, even when $\lambda$ is not a
unit modulo $N$.
\end{lemma}

\begin{proof}
For prime $\lambda$, this is \cite[Lemma~3.1]{CovenMeyerowitz}.
The tiling identity gives $|\lambda U||V|=N=|U||V|$, so multiplication
by $\lambda$ is injective on $U$. Applying the prime case successively
to the prime factors of an arbitrary $\lambda$ proves the result,
since the cardinality of the first factor remains $|U|$ at each step.
\end{proof}

\begin{lemma}\label{lem:forbidden}
Suppose $A\oplus B=\mathbb{Z}_{p^e m}$, where $p\nmid m$.
Let $\delta$ be a nonzero difference from $A$ and $\beta$ a nonzero
difference of integer representatives of $B$.
It is impossible to have both
\begin{equation}\label{eq:forbidden}
v_p(\delta)=v_p(\beta)<e,
\qquad \gcd(m,\delta)\mid\beta.
\end{equation}
\end{lemma}

\begin{proof}
Write $s=v_p(\delta)$ and $g=\gcd(m,\delta)$.
The conditions imply $\gcd(\delta,p^e m)=p^s g\mid\beta$, so
$\lambda\delta\equiv\beta\pmod{p^e m}$ has a positive integer solution $\lambda$.
Dividing this congruence by $p^s$ and reducing modulo $p$ gives
$p\nmid\lambda$. It creates a collision in $\lambda A+B$, contradicting
\eqref{eq:dilation}. The difference $\beta$ is nonzero modulo $p^e$,
so the two terms from each factor are distinct.
\end{proof}

The next lemma is the finite cyclic form of complementary branching
\cite[Theorems~4.2 and~5.1(b)]{FanFanShi}. We give a direct proof from
Lemma~\ref{lem:forbidden}.

\begin{lemma}\label{lem:pure}
Let $e\geq n$. The tilings of $A$ in $\mathbb{Z}_{p^e}$ are
exactly the nonempty sets $C\subset\mathbb{Z}_{p^e}$ such that, for
each $1\leq j\leq e$, reduction from $r_{p^j}(C)$ onto
$r_{p^{j-1}}(C)$ is bijective when $j\in I$, and has every fibre of
size $p$ when $j\notin I$. Equivalently, $C$ has single branching
at the positions in $I$ and full branching at all other positions.
\end{lemma}

\begin{proof}
Suppose $A\oplus C=\mathbb{Z}_{p^e}$. If an occupied class of $C$
has two lifts at position $\nu_i$, their extensions to $C$ give a
difference of valuation $\nu_i-1$. A difference of $A$ with that
valuation exists by \eqref{eq:levels}, contradicting
Lemma~\ref{lem:forbidden} with $m=1$. Thus $C$ has single branching
at each position in $I$. At the other $e-k$ positions there are at
most $p$ lifts per occupied class. Equality in $|C|=p^{e-k}$ forces
full branching at every such class.

Conversely, the stated branching gives $|C|=p^{e-k}$.
Differences of distinct integer representatives of $C$ have valuations
$j-1$ with $j\notin I$, all less than $e$. These valuations are disjoint
from those of differences from $A$, so no collision can occur modulo
$p^e$. Hence the addition map $A\times C\to\mathbb{Z}_{p^e}$ is injective,
and the cardinalities make it bijective.
\end{proof}

For such a $C$, call the positions in $I$ selected and the other
positions free. A selected digit can depend on the lower free digits.
For example, in $\{0,2\}\oplus\{0,3\}=\mathbb{Z}_4$, the upper binary
digit of the second factor equals its lower digit.

For the gcd calculation, define
\[
\mathcal{E}_i=\{a-a':a,a'\in A,\ a\neq a',\ v_p(a-a')=\nu_i-1\}.
\]

\begin{lemma}\label{lem:gcd}
The integers in \eqref{eq:gcd-data} satisfy \eqref{eq:chain}.
Furthermore,
\begin{equation}\label{eq:exact-gcd}
\gcd\{|\delta|:\delta\in\mathcal{E}_i\}=Q_i\Delta_i.
\end{equation}
\end{lemma}

\begin{proof}
A pair that first differs at position $\nu_i$ supplies a normalized
difference prime to $p$. Thus $\Delta_i$ is defined and $p\nmid\Delta_i$.
Each difference used at level $i+1$ is also used at level $i$ before
normalization. Hence
$\Delta_i\mid p^{\nu_{i+1}-\nu_i}\Delta_{i+1}$, proving the chain.

To prove \eqref{eq:exact-gcd}, consider two elements in one class
modulo $Q_i$. If they also agree modulo $pQ_i$, choose an element
$a''$ in another one of the $p$ branches of that class.
Their difference is the sum of two differences through $a''$, both
of exact valuation $\nu_i-1$. Thus all differences in the definition
of $Q_i\Delta_i$ are generated by $\mathcal{E}_i$, and the two gcds
are equal.
\end{proof}

\section{Possible fundamental periods}\label{sec:periods}

For a residue set $B\subset\mathbb{Z}_N$, form $B(x)$ using its
representatives in $\{0,\ldots,N-1\}$. With $G_N(x)=\sum_{j=0}^{N-1}x^j$,
the tiling identity in $\mathbb{Z}[x]/(x^N-1)$ is
\begin{equation}\label{eq:group-ring}
A\oplus B=\mathbb{Z}_N
\quad\Longleftrightarrow\quad
A(x)B(x)\equiv G_N(x)\pmod{x^N-1}
\end{equation}
\cite[Lemma~1.3]{CovenMeyerowitz}.

We first identify the cyclotomic divisors of $A(x)$ whose orders are
multiples of $P$.

\begin{proposition}\label{prop:cyclotomic}
If $\alpha\geq n$, $m\geq1$, and $p\nmid m$, then
\begin{equation}\label{eq:cyclotomic}
\Phi_{p^\alpha m}(x)\mid A(x)
\quad\Longleftrightarrow\quad
\alpha=n\ \text{and}\ m\mid D.
\end{equation}
\end{proposition}

\begin{proof}
Let $\xi,\eta$ have orders $p^\alpha,m$. Their product has order
$p^\alpha m$, so the divisibility in \eqref{eq:cyclotomic} is equivalent
to $A(\xi\eta)=0$.
The minimal polynomial of $\xi$ over $\mathbb{Q}(\eta)$ is
$\Phi_{p^\alpha}$, since the coprime conductors give extension
degree $\varphi(p^\alpha)$, where $\varphi$ is Euler's totient function.
Reduce $A(\xi\eta)$ by collecting exponents modulo $p^\alpha$.
The resulting polynomial has degree less than $p^\alpha$, with
each coefficient either zero or a single root $\eta^a$, because
$A$ is injective modulo $P$.
The quotient by $\Phi_{p^\alpha}$, when it exists, has degree less
than $p^{\alpha-1}$. Thus the polynomial vanishes at $\xi$ exactly
when its coefficients agree in every block
\[
r,\ r+p^{\alpha-1},\ \ldots,\ r+(p-1)p^{\alpha-1},
\qquad 0\leq r<p^{\alpha-1}.
\]
For $\alpha>n$, a nonzero block would give distinct elements of
$A$ congruent modulo $P$, which is impossible.
For $\alpha=n$, each block indexed by an occupied class modulo $Q_k$
has all $p$ coefficients nonzero. Its coefficients agree precisely when all differences within
each class modulo $Q_k$ are divisible by $m$.
By \eqref{eq:gcd-data} and $p\nmid m$, this means $m\mid D$.
\end{proof}

To prove Theorem~\ref{thm:periods}, we compare the mask of a tiling
with the mask of its translate. A related use of cyclotomic quotients
appears in \cite[Lemma~2]{LabaZakharov}. We use the standard description
of the prime-power cyclotomic orders
\[
S_A=\{q^j:q\text{ prime},\ j\geq1,\ \Phi_{q^j}(x)\mid A(x)\}
=\{p^{\nu_1},\ldots,p^{\nu_k}\}
\]
\cite[Theorem~B1 and Section~4]{CovenMeyerowitz}.

\begin{proof}[Proof of Theorem~\ref{thm:periods}]
Let $T$ be a tiling with a period $N=p^e m$, where $p\nmid m$.
Every period is divisible by $\operatorname{lcm}(S_A)=P$
\cite[remarks following Lemma~2.1]{CovenMeyerowitz}, so $e\geq n$.
Write $B=r_N(T)$ and $L=Q_km$. Let $F(x)$ be the representative of
$(x^{PD}-1)B(x)$ modulo $x^N-1$ of degree less than $N$.
For $s\mid N$ with $s\nmid L$, the exponent of $p$ in $s$ is at least $n$.
By \eqref{eq:group-ring}, $\Phi_s(x)$ divides $A(x)$ or $B(x)$.
In the first case Proposition~\ref{prop:cyclotomic} gives $s\mid PD$,
so in either case $\Phi_s(x)$ divides $F(x)$. Therefore
\[
\frac{x^N-1}{x^L-1}=\sum_{j=0}^{N/L-1}x^{jL}\mid F(x).
\]
The quotient has degree less than $L$, so the coefficients of $F(x)$
are constant on residue classes modulo $L$. These coefficients are the
differences of the indicators of $B+PD$ and $B$. A nonzero coefficient
would give two standard representatives differing by $L$ in one of these two sets.
Any $\delta\in\mathcal{E}_k$ has $v_p(\delta)=v_p(L)=n-1<e$ and
$\gcd(m,\delta)\mid L$, so Lemma~\ref{lem:forbidden} rules out that pair
in either tiling. Hence $F(x)=0$ and $PD$ is a period of $T$.
The fundamental period is a divisor of $PD$ and a multiple of $P$.

For the converse, fix $d\mid D$ and define
\begin{equation}\label{eq:short-realization}
B=\left\{\sum_{j\in J}t_jp^{j-1}+Q_k z:
t_j\in\{0,\ldots,p-1\}\ (j\in J),\ z\in\{0,\ldots,d-1\}\right\}
\subset\mathbb{Z}_{Pd}.
\end{equation}
The sums over $J$ lie in $[0,Q_k)$ since $n\in I$. In a collision
$a+b\equiv a'+b'\pmod{Pd}$, unequal sums over $J$ would force the first
differing digit of $b,b'$ to lie in $J$. The congruence modulo $P$ would
then give the same first differing digit for $a,a'$, contrary to
\eqref{eq:levels}. Thus the two sums over $J$ agree and $a\equiv a'\pmod{Q_k}$.
By \eqref{eq:gcd-data}, $Q_kD$ divides $a-a'$. The congruence modulo $d$
now gives $z=z'$ since $\gcd(Q_k,d)=1$ and $0\leq z,z'<d$.
Then $a=a'$ by injectivity modulo $P$. The cardinalities give
$|A||B|=p^k p^{n-k}d=Pd$, so $B$ tiles with $A$ modulo $Pd$.
Moreover, $B$ contains $0$ and lies in $[0,Pd/p)$. If its fundamental
period were a proper divisor $h$ of $Pd$, then $h\leq Pd/2$ and
$Pd-h$ would belong to $B$, which is impossible. Thus every $Pd$ occurs.
The cyclotomic description in \eqref{eq:period-set} follows from
Proposition~\ref{prop:cyclotomic}.
\end{proof}

\section{Compatibility of the CRT fibres}\label{sec:classification}

By Theorem~\ref{thm:periods}, every integer tiling is the extension of a
tiling at a modulus $Pd$ with $d\mid D$. We now prove the criterion at
each of these moduli.

\begin{proof}[Proof of Theorem~\ref{thm:classification}]
Suppose first that $A\oplus B=\mathbb{Z}_{Pd}$.
Choose a positive integer $\lambda$ with $\lambda\equiv1\pmod P$ and
$\lambda\equiv0\pmod d$. Dilation by $\lambda$ gives
\eqref{eq:fibre-tiling}, since $\rho_d(\lambda A)=\overline{A}\times\{0\}$.

We prove \eqref{eq:compatibility} by induction on $i$.
Fix $\delta\in\mathcal{E}_i$ and define $g=\gcd(d,\delta)$.
If $i>1$, then $d_{i-1}\mid g$ because $\Delta_{i-1}\mid\Delta_i\mid\delta$.
The induction hypothesis gives equality of the projections of $C_y$
and $C_{y+g}$ at position $\nu_{i-1}$. All positions between $\nu_{i-1}$
and $\nu_i$ have full branching by Lemma~\ref{lem:pure}. Hence the two
fibres have the same image modulo $Q_i$. For $i=1$ this image is all
of $\mathbb{Z}_{Q_1}$.

If $\pi_i(C_y)\neq\pi_i(C_{y+g})$, their unique lifts of some common
class modulo $Q_i$ differ modulo $pQ_i$. Choose points of $B$
extending these lifts. Their difference $\beta$, taken in integer
representatives, satisfies
$v_p(\beta)=\nu_i-1$ and $g\mid\beta$.
The latter follows by reduction modulo $d$, where the second
coordinates differ by $g$. This contradicts
Lemma~\ref{lem:forbidden}. Thus $\pi_i(C_y)$ is invariant under every
shift $\gcd(d,\delta)$ with $\delta\in\mathcal{E}_i$.
These shifts generate the subgroup with positive generator
\[
\gcd\bigl(d,\gcd\{|\delta|:\delta\in\mathcal{E}_i\}\bigr)
=\gcd(d,Q_i\Delta_i)=d_i,
\]
by Lemma~\ref{lem:gcd}. This proves the required synchronization.

Conversely, suppose \eqref{eq:fibre-tiling} and
\eqref{eq:compatibility} hold. A collision $a+b=a'+b'$ modulo $Pd$
with $a\neq a'$ has $v_p(a-a')=\nu_i-1$ for some $i$.
Write $\rho_d(b)=(x,y)$ and $\rho_d(b')=(x',z)$.
Since $\Delta_i\mid a-a'$ and $d_i\mid d$, one has
$y\equiv z\pmod{d_i}$. Also $x\equiv x'\pmod{Q_i}$.
By \eqref{eq:compatibility}, their images modulo $pQ_i$ belong to the
same projected fibre, and reduction of that fibre modulo $Q_i$ is
injective by Lemma~\ref{lem:pure}. Thus
$x\equiv x'\pmod{pQ_i}$, contradicting $v_p(a-a')=\nu_i-1$.
Finally, \eqref{eq:fibre-tiling} gives $|B|=p^{n-k}d$, so the sum
map is bijective.

Within a constant run of the $d_i$, equality at the last projection
implies equality at the preceding projections. The last run has value
$d_k=d$, for which $y\equiv z\pmod{d}$ means $y=z$.

Every period of $B+Pd\mathbb{Z}$ is a multiple of $P$ by
Theorem~\ref{thm:periods}. Under $\rho_d$, a shift by $Ph$ fixes the first
coordinate and translates the second by $r_d(Ph)$. Thus the admissible
shifts by multiples of $P$ are exactly those for which $r_d(Ph)\in H_B$.
The inverse image of $H_B$ in $\mathbb{Z}$ is $(d/|H_B|)\mathbb{Z}$,
and $P$ is coprime to $d/|H_B|$. This proves
\eqref{eq:fibre-fundamental-period}.
\end{proof}

For example, let $A=\{0,3,10,13\}$. Its data are
\[
p=2,\qquad(\nu_1,\nu_2)=(1,2),\qquad
(\Delta_1,\Delta_2)=(1,5),\qquad P=4,\quad D=5.
\]
Since $\overline{A}=\mathbb{Z}_4$, every $C_y$ is a singleton.
The only compatibility condition requires these singletons to have
the same parity. Thus all tilings have the form
\begin{equation}\label{eq:running-example}
\rho_5(B)=\{(c+2g(y),y):y\in\mathbb{Z}_5\},
\qquad c\in\{0,1\},\quad g:\mathbb{Z}_5\longrightarrow\{0,1\}.
\end{equation}
There are $2\cdot2^5=64$ choices. The four with constant $g$ have
fundamental period $4$; the other sixty have fundamental period $20$.

\section{Digit parameters and counting}\label{sec:counting}

Write $\mathcal{C}_A(N)=\{B\subset\mathbb{Z}_N:A\oplus B=\mathbb{Z}_N\}$.
Let $\mathcal{D}=\{0,\ldots,p-1\}$ and
$J_i=\{j\in J:j<\nu_i\}$, so $|J_i|=\nu_i-i$.
For a finite index set $S$, $\mathcal{D}^S$ denotes the set of
digit tuples indexed by $S$, including the empty tuple when $S$ is empty.
Write $\mathbf{t}|_S$ for restriction to $S$.

Fix $d\mid D$. At position $\nu_i$, Lemma~\ref{lem:pure} gives a unique
digit above each occupied residue class modulo $Q_i$ in a fibre $C_y$.
We call the lower digits specifying such a class a prefix.
Theorem~\ref{thm:classification} identifies these digits across fibres
whose second coordinates agree modulo $d_i$. The lower free digits
index the occupied prefixes, giving the following parametrization.

\begin{proposition}\label{prop:digit-classification}
For $1\leq i\leq k$, choose a function
\begin{equation}\label{eq:functions}
f_i:\mathcal{D}^{J_i}\times\mathbb{Z}_{d_i}\longrightarrow\mathcal{D}.
\end{equation}
For $f=(f_1,\ldots,f_k)$, $\mathbf{t}\in\mathcal{D}^J$ and
$y\in\mathbb{Z}_d$, define
\begin{align}
x_f(\mathbf{t},y)&=\sum_{j\in J}t_jp^{j-1}
+\sum_{i=1}^k f_i(\mathbf{t}|_{J_i},r_{d_i}(y))p^{\nu_i-1},
\label{eq:digit-form}\\
B_f&=\rho_d^{-1}\{(x_f(\mathbf{t},y),y):\mathbf{t}\in\mathcal{D}^J,
\ y\in\mathbb{Z}_d\}.
\label{eq:digit-set}
\end{align}
The map $f\mapsto B_f$ is a bijection onto $\mathcal{C}_A(Pd)$.
\end{proposition}

\begin{proof}
In each fibre, full branching at free positions and single branching
at selected positions identify the occupied prefixes modulo $Q_i$
bijectively with $\mathcal{D}^{J_i}$. For a given lower free tuple,
the selected digit at $\nu_i$ is determined by its occupied prefix
and the fibre coordinate $y$. If $y\equiv z\pmod{d_i}$,
\eqref{eq:compatibility} gives equality of the two projections modulo
$p^{\nu_i}$. The two prefixes labelled by the same lower free tuple
therefore have the same selected digit. This defines $f_i$ on the
domain in \eqref{eq:functions}, and reading the successive digits gives
\eqref{eq:digit-form}.

Conversely, arbitrary functions on these domains give full branching
at the free positions and single branching at the selected positions,
so \eqref{eq:fibre-tiling} follows from Lemma~\ref{lem:pure}.
If $y\equiv z\pmod{d_i}$, then $y\equiv z\pmod{d_h}$ for every $h\leq i$.
Thus for each free tuple, the digits through position $\nu_i$ agree in
the two fibres. This gives \eqref{eq:compatibility}. Every lower free
tuple and residue modulo $d_i$ occurs because restriction of tuples
and reduction from $\mathbb{Z}_d$ are surjective. Hence $B_f$ determines
every value of every $f_i$.
\end{proof}

\begin{example}\label{ex:free-digit}
For $A=\{0,3,20,23\}$, the levels are $(\nu_1,\nu_2)=(1,3)$
and the gcd chain is $(1,5)$. Thus $P=8$, $D=5$, and $J=\{2\}$.
In $\mathbb{Z}_8\times\mathbb{Z}_5$, all tilings are
\[
\{(c+2t+4f(t,y),y):t\in\{0,1\},\ y\in\mathbb{Z}_5\},
\]
where $c\in\{0,1\}$ and $f:\{0,1\}\times\mathbb{Z}_5\longrightarrow\{0,1\}$
are arbitrary. The chosen upper digit can depend on both the lower
free digit and the coprime coordinate.
\end{example}

We count tilings with a given period and tilings with that fundamental
period separately. For $N\geq1$, define
\[
F_A(N)=|\mathcal{C}_A(N)|,\qquad
\mathcal{E}_A(N)=\{B\in\mathcal{C}_A(N):\operatorname{per}(B+N\mathbb{Z})=N\},
\qquad E_A(N)=|\mathcal{E}_A(N)|.
\]
Write $\mu$ for the M\"obius function.

\begin{corollary}\label{cor:counts}
For every positive integer $N$,
\begin{equation}\label{eq:all-counts}
F_A(N)=
\begin{cases}
0,&P\nmid N,\\[2pt]
p^{\sum_{i=1}^k p^{\nu_i-i}\gcd(N,\Delta_i)},&P\mid N.
\end{cases}
\end{equation}
For $d\mid D$,
\begin{equation}\label{eq:exact-counts}
E_A(Pd)=\sum_{e\mid d}\mu(d/e)F_A(Pe).
\end{equation}
The total number of integer tilings is
\begin{equation}\label{eq:total-count}
F_A(PD)=p^{\sum_{i=1}^k p^{\nu_i-i}\Delta_i}.
\end{equation}
There are $E_A(Pd)/(Pd)$ translation classes of fundamental period $Pd$.
The total number of translation classes is
\begin{equation}\label{eq:total-translations}
\frac{1}{PD}\sum_{d\mid D}\varphi(D/d)F_A(Pd).
\end{equation}
\end{corollary}

\begin{proof}
By Theorem~\ref{thm:periods}, there are no tilings if $P\nmid N$.
For $P\mid N$, let $d=\gcd(N,D)$. Every tiling with period $N$ has
period $\gcd(N,PD)=Pd$, and every tiling with period $Pd$ has period $N$.
This gives a bijection of the residue sets by reduction from $N$ to $Pd$
and full inverse image in the other direction.
Proposition~\ref{prop:digit-classification} gives $p^{p^{\nu_i-i}d_i}$
choices for $f_i$, since its domain has $p^{\nu_i-i}d_i$ elements.
Multiplication gives \eqref{eq:all-counts}, using
$d_i=\gcd(d,\Delta_i)=\gcd(N,\Delta_i)$.
Every tiling counted by $F_A(Pd)$ has fundamental period $Pe$ for
exactly one $e\mid d$. M\"obius inversion gives
\eqref{eq:exact-counts}.
The common period $PD$ identifies all integer tilings with
their unique residue sets at that modulus, proving
\eqref{eq:total-count}. Translation acts freely on the stratum of fixed fundamental period.
Summing $E_A(Pd)/(Pd)$ over $d\mid D$, substituting
\eqref{eq:exact-counts}, and using
$\sum_{h\mid q}\mu(h)/h=\varphi(q)/q$ gives
\eqref{eq:total-translations}.
\end{proof}

A tiling $B\in\mathcal{C}_A(Pd)$, with $d\mid D$, has fundamental period $Pd$ if and only if
\begin{equation}\label{eq:period-test}
B+Pd/q\neq B\quad\text{for every prime }q\mid d.
\end{equation}
Indeed, if the fundamental period is $Pe$ with $e\mid d$ and $e<d$,
choose a prime $q\mid d/e$. Then $Pe\mid Pd/q$, so $Pd/q$ is a period.
Conversely, invariance under $Pd/q$ gives a period smaller than $Pd$.
For $d=1$ the condition is empty and every
tiling has fundamental period $P$.

For an individual $B_f$, the digit functions also give a way to evaluate
\eqref{eq:fibre-fundamental-period}. For each $i$, let $e_i$ be the
fundamental period in the second variable of $f_i$, with the same shift
required to preserve its values for all lower free tuples. The shifts
preserving $C_y$ for every $y$ are precisely the shifts preserving all
the functions $f_i$. Thus the inverse image of $H_{B_f}$ in $\mathbb{Z}$
is $\operatorname{lcm}(e_1,\ldots,e_k)\mathbb{Z}$, and
\[
\operatorname{per}(B_f+Pd\mathbb{Z})=P\operatorname{lcm}(e_1,\ldots,e_k).
\]

\begin{example}\label{ex:lcm}
Let
\[
A=\{0,3,30,33\}=\{0,3\}\oplus\{0,30\}.
\]
Here $P=4$, $(\nu_1,\nu_2)=(1,2)$, and
$(\Delta_1,\Delta_2)=(3,15)$. Every tiling at modulus $60$
is a graph
\[
\{(f_1(y\bmod3)+2f_2(y),y):y\in\mathbb{Z}_{15}\},
\]
where $f_1:\mathbb{Z}_3\longrightarrow\{0,1\}$ and
$f_2:\mathbb{Z}_{15}\longrightarrow\{0,1\}$ are arbitrary.
Choose $f_1$ to be the indicator of $0\bmod3$, and $f_2$ the
indicator of $0\bmod5$. Their fundamental periods are $3$ and $5$,
so the tiling has fundamental period $4\operatorname{lcm}(3,5)=60$.
Neither function has fundamental period $15$.
\end{example}

\section{Partitions into translates}\label{sec:partitions}

Fix $d\mid D$ and $M=Pd$. A tiling $B\in\mathcal{C}_A(M)$ induces the
unordered partition
\[
\mathcal{P}_M(B)=\{A_M+b:b\in B\},\qquad A_M=r_M(A),
\]
of $\mathbb{Z}_M$. We call $B$ and $B'$ \emph{equivalent as tile
partitions} if $\mathcal{P}_M(B)=\mathcal{P}_M(B')$. We also consider
orbits of these partitions under simultaneous translation of their tiles.
The fundamental period of a partition is the least positive integer
whose translation preserves the partition. It divides $M$ and can be
smaller than the fundamental period of the tiling that induces it.

\begin{theorem}\label{thm:tile-partitions}
Define
\begin{equation}\label{eq:tile-stabilizer-order}
s=\max\{j\in\{1,\ldots,k\}:
 \nu_{k-j+1}=n-j+1,\ d\mid\Delta_{k-j+1}\},
\qquad r=k-s,\qquad L=M/p^s .
\end{equation}
The translation stabilizer of $A_M$ is $L\mathbb{Z}_M$. For tilings
$B,B'\in\mathcal{C}_A(M)$,
\[
\mathcal{P}_M(B)=\mathcal{P}_M(B')
\quad\Longleftrightarrow\quad r_L(B)=r_L(B').
\]
Every such partition is induced by a tiling in $\mathcal{E}_A(M)$.
The distinct partitions are parametrized by the first $r$ functions
in \eqref{eq:functions}, chosen independently. To obtain $r_L(B)$,
use \eqref{eq:digit-form} with the sum of selected digits truncated at
$i=r$. Apply the inverse Chinese remainder isomorphism to the resulting
pairs in $\mathbb{Z}_{p^{n-s}}\times\mathbb{Z}_d$.
The number of partitions is
\begin{equation}\label{eq:tile-partition-total}
p^{\sum_{i=1}^{r}p^{\nu_i-i}d_i}.
\end{equation}
In particular there is exactly one partition if and only if $r=0$.
\end{theorem}

\begin{proof}
The order of the translation stabilizer divides $|A_M|=p^k$. In CRT
coordinates, the subgroup of order $p^j$ varies the highest $j$ digits
while fixing both the residue modulo $p^{n-j}$ and the coordinate modulo
$d$. Reduction of $A$ modulo $P$ is injective. By the branching description
in Section~\ref{sec:digits}, invariance under this subgroup requires all
highest $j$ positions to belong to $I$. This is equivalent to
$\nu_{k-j+1}=n-j+1$. With these positions in $I$, all $p^j$ lifts of an
occupied prefix modulo $p^{n-j}$ occur. They lie in one coset of the
subgroup precisely when the corresponding elements of $A$ agree modulo
$d$. Since $Q_{k-j+1}=p^{n-j}$ is coprime to $d$, this condition on all
occupied prefixes is equivalent to $d\mid\Delta_{k-j+1}$. This proves the
stabilizer formula. The index $j=1$ always satisfies these conditions,
so $s\geq1$.

The equality $A_M+b=A_M+b'$ holds exactly when $b-b'$ belongs to the
stabilizer. Thus the image $r_L(B)$ determines the partition, and the
partition recovers that image. Reduction modulo $L$ erases the highest
$s$ digits in the first CRT coordinate and preserves the second
coordinate. The erased digits are all selected. Within each fibre,
the lower free tuples still label each occupied prefix uniquely, and the surviving digits
recover every value of the first $r$ functions. Conversely, any choice
of these functions extends to a tiling at modulus $M$ by
Proposition~\ref{prop:digit-classification}.

Let $\widetilde B$ consist of the representatives of $r_L(B)$ in $[0,L)$.
Reduction modulo $L$ is injective on $B$, since two elements in the
same class would give the same tile. Each replacement changes a
translate by an element of the stabilizer, so $\widetilde B$ induces
the same partition. If its periodic extension had a proper fundamental
period $h\mid M$, then $h\leq M/2$. Since $L\leq M/2$, translating
the largest element of $\widetilde B$ by $h$ would give a larger
element below $M$, a contradiction. Thus $\widetilde B\in\mathcal{E}_A(M)$.
Finally the domain of $f_i$ has $p^{\nu_i-i}d_i$ elements. Multiplying
the independent choices gives \eqref{eq:tile-partition-total}.
\end{proof}

If $r>0$, write $P_0=p^{\nu_r}$ and $d_*=d_r$. If $r=0$, set
$P_0=d_*=1$. The quantities below are for the fixed modulus $M$, so
the number $r$ of retained functions stays fixed throughout.

\begin{corollary}\label{cor:tile-partition-counts}
For $e\mid d_*$, the number of partitions with period $P_0e$ and the
number with fundamental period $P_0e$ are, respectively,
\begin{align}
F(e)&=p^{\sum_{i=1}^{r}p^{\nu_i-i}\gcd(e,d_i)},\label{eq:tile-partition-F}\\
E(e)&=\sum_{a\mid e}\mu(e/a)F(a).\label{eq:tile-partition-E}
\end{align}
The possible fundamental periods are exactly the $P_0e$ with $e\mid d_*$.
The total number of
translation classes of partitions is
\begin{equation}\label{eq:tile-partition-orbits}
\sum_{e\mid d_*}\frac{E(e)}{P_0e}
=\frac{1}{P_0d_*}\sum_{e\mid d_*}\varphi(d_*/e)F(e).
\end{equation}
For one fixed partition of fundamental period $P_0e$, the number of
tilings of fundamental period $M$ that induce it is
\begin{equation}\label{eq:tile-partition-lifts}
\Lambda(e)=\sum_{\substack{a\mid d\\ e\mid a}}
             \mu(d/a)\,p^{s\,p^{n-k}a}.
\end{equation}
There are $\Lambda(e)P_0e/M$ translation classes of tilings above its
translation class of partitions.
\end{corollary}

\begin{proof}
A shift preserves the partition exactly when it preserves $r_L(B)$.
If $r=0$, the projected set is the whole quotient. If $r>0$, every
projected digit position above $\nu_r$ is free and the retained
functions depend on $y$ only modulo $d_*$. Thus $r_L(B)$ is the full
inverse image of a set modulo $P_0d_*$. A nontrivial $p$-part of its
stabilizer at that modulus would contain the subgroup of order $p$.
This subgroup fixes $y$ and the prefix below $\nu_r$, and varies
the digit at $\nu_r$, where the lift is unique. Therefore each
period is a multiple of $P_0$.

For $e\mid d_*$, a shift by $P_0e$ fixes the first coordinate modulo
$P_0$. It preserves the projected set exactly when all retained
functions are invariant under the shift by $P_0e$ in the second
variable. Since $P_0$ is coprime to $d_i$, this means that $f_i$ factors
through reduction modulo $\gcd(e,d_i)$. Counting these functions gives
\eqref{eq:tile-partition-F}. When $r>0$, constant functions for $i<r$ and
the indicator of $y\equiv0\pmod e$ for $f_r$ realize fundamental period
$P_0e$. When $r=0$ the sole partition has period $1$. M\"obius inversion gives
\eqref{eq:tile-partition-E}, and summing the orbit counts as in
\eqref{eq:total-translations} gives \eqref{eq:tile-partition-orbits}.

Fix a partition of fundamental period $P_0e$. Since $P_0\mid P$,
its retained functions can occur in a tiling of period $Pa$, for $a\mid d$,
exactly when $e\mid a$. For each of the last $s$ indices we have
$d_i=d$ and $J_i=J$, with $|J|=n-k$. At period $Pa$, each of these
functions has $p^{n-k}a$ independent entries. Hence there are
$p^{s\,p^{n-k}a}$ choices when $e\mid a$. Inverting over $a\mid d$
gives \eqref{eq:tile-partition-lifts}. The subgroup of translations
preserving the fixed partition has order $M/(P_0e)$ and acts freely
on its exact-period-$M$ tilings. Its orbits are precisely the
translation classes of tilings above the partition orbit.
\end{proof}

\begin{figure}[!htbp]
	\centering
	\begin{tikzpicture}
		\definecolor{partitionTile1}{HTML}{CBD4DF}
		\definecolor{partitionTile2}{HTML}{589FC8}
		\definecolor{partitionTile3}{HTML}{B9D3BA}
		\definecolor{partitionTile4}{HTML}{D5CBE4}
		\definecolor{partitionTile5}{HTML}{F6C549}
		\definecolor{partitionTile6}{HTML}{E6CBBE}
		\definecolor{partitionInner1}{HTML}{B3C4D3}
		\definecolor{partitionInner2}{HTML}{E1E5E9}
		\definecolor{partitionInk}{HTML}{232D3A}
		
		\foreach \InnerRadius/\OuterRadius/\RingModulus/\RingStarts in {
			2.69/3.64/24/{0,1,2,6,7,8},
			1.51/2.45/24/{0,10,2,6,4,8},
			0.39/1.27/8/{0,2}} {
			\pgfmathsetmacro{\LabelRadius}{(\InnerRadius+\OuterRadius)/2}
			\foreach \Start [count=\TileNumber from 1] in \RingStarts {
				\edef\TileFill{partitionTile\TileNumber}
				\ifnum\RingModulus=8\relax
				\edef\TileFill{partitionInner\TileNumber}
				\fi
				\foreach \TileOffset in {0,3,12,15} {
					\pgfmathtruncatemacro{\Residue}{mod(\Start+\TileOffset,\RingModulus)}
					\pgfmathsetmacro{\SectorStart}{90-360*\Residue/\RingModulus}
					\pgfmathsetmacro{\SectorEnd}{90-360*(\Residue+1)/\RingModulus}
					\pgfmathsetmacro{\SectorMiddle}{(\SectorStart+\SectorEnd)/2}
					\draw[fill=\TileFill,draw=white,line width=0.4pt]
					(\SectorStart:\InnerRadius) -- (\SectorStart:\OuterRadius)
					arc[start angle=\SectorStart,end angle=\SectorEnd,radius=\OuterRadius cm]
					-- (\SectorEnd:\InnerRadius)
					arc[start angle=\SectorEnd,end angle=\SectorStart,radius=\InnerRadius cm]
					-- cycle;
					\node[font=\scriptsize\sffamily\bfseries,text=partitionInk,inner sep=0pt]
					at (\SectorMiddle:\LabelRadius) {\TileNumber};
				}
			}
			\draw[partitionInk,line width=0.4pt] (0,0) circle[radius=\InnerRadius cm];
			\draw[partitionInk,line width=0.4pt] (0,0) circle[radius=\OuterRadius cm];
		}
		\foreach \Residue in {4,10,16,22} {
			\pgfmathsetmacro{\MarkerAngle}{90-15*(\Residue+0.5)}
			\fill[partitionInk] (\MarkerAngle:2.57) circle[radius=1.6pt];
		}
		\node[font=\scriptsize,text=partitionInk] at (82.5:3.91) {$0$};
		\begin{scope}[
			every node/.style={anchor=west,align=left,inner sep=0pt,
				font=\small,text=partitionInk}]
			\node at (4.05,1.85) {\textbf{Outer ring} $\boldsymbol{(B_1)}$\\
				Fundamental tiling period: $24$\\Fundamental partition period: $6$};
			\node at (4.05,0) {\textbf{Middle ring} $\boldsymbol{(B_2)}$\\
				Fundamental tiling period: $24$\\Fundamental partition period: $2$};
			\node at (4.05,-1.85) {\textbf{Inner ring} $\boldsymbol{(B_3)}$\\
				Fundamental tiling period: $8$\\Fundamental partition period: $2$};
		\end{scope}
	\end{tikzpicture}
	\caption{Representatives of the translation classes of tile partitions
		in Example~\ref{ex:tile-partitions-24}, grouped by fundamental tiling
		period. In each ring the residues increase clockwise from the top
		divider, starting at $0$. Wedges with the same number in a ring form
		one translate of $A$. In the two $24$-sector rings, the dots mark the
		residues $4,10,16,22$ that switch between tiles $2$ and $5$.}
	\label{fig:tile-partitions-rings}
\end{figure}
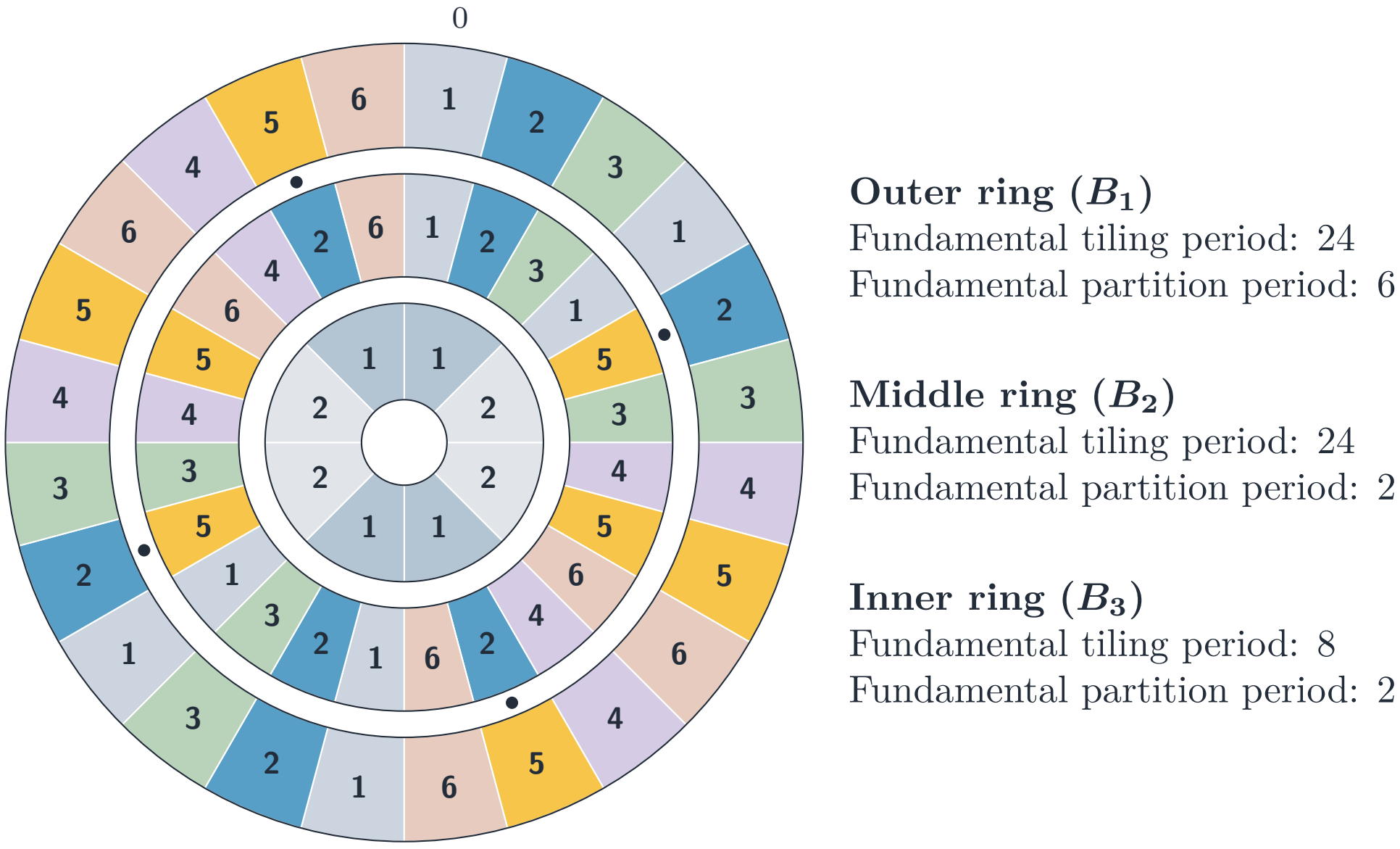

\begin{example}\label{ex:tile-partitions-24}
	Let $A=\{0,3,12,15\}$. Here $P=8$, $D=3$, and
	$(\nu_1,\nu_2)=(1,3)$, so the possible fundamental tiling periods
	are $8$ and $24$. In both cases $s=r=1$ and $P_0=2$.
	
	At $M=24$ we have $L=12$ and $d_*=3$. The two partitions of
	fundamental period $2$ form one translation class, and the six of
	fundamental period $6$ form another. Representatives of these classes
	with fundamental tiling period $24$ are
	\[
	B_1=\{0,1,2,6,7,8\},\qquad B_2=\{0,2,4,6,8,10\}.
	\]
	Modulo $12$, we have $B_i\oplus\{0,3\}=\mathbb{Z}_{12}$.
	Since $r_{24}(A)=\{0,3\}\oplus\{0,12\}$, it follows that
	$A\oplus B_i=\mathbb{Z}_{24}$ for $i=1,2$. Both sets lie below $12$
	and hence have fundamental period $24$ by the argument in
	Theorem~\ref{thm:tile-partitions}. Their projected sets have
	fundamental periods $6$ and $2$, respectively. Each fixed partition
	of period $2$ has $60$ exact-period-$24$ tilings, and each of period
	$6$ has $64$.
	
	At $M=8$ we have $L=4$ and $d_*=1$. The two parity classes of starts
	modulo $4$ give two partitions of fundamental period $2$ in one
	translation class. Each is induced by four tilings of fundamental
	period $8$. One representative is $B_3=\{0,2\}$. Its two translates
	of $A$ modulo $8$ partition $\mathbb{Z}_8$, and $B_3\subset[0,4)$
	has fundamental period $8$. Figure~\ref{fig:tile-partitions-rings}
	displays one representative of each class at its specified
	fundamental tiling period. The distinct images of the middle ring's
	tiles under reduction modulo $8$ form the inner partition.
\end{example}

For the running example $A=\{0,3,10,13\}$ at $M=20$, we have
$L=10$, $P_0=2$ and $d_*=1$. The projected sets are the even and odd
residues modulo $10$. They give two partitions in one translation
class, each induced by $30$ tilings of fundamental period $20$.

\section{The canonical representative and sharp bounds}\label{sec:canonical}

\begin{proof}[Proof of Theorem~\ref{thm:canonical}]
Set $g_i=Q_i\Delta_i$. Then $g_i\mid g_{i+1}$ and
$v_p(g_i)=\nu_i-1$.
In the difference of two sums in \eqref{eq:canonical}, let $i$
be the first index with a nonzero digit difference.
That coefficient is nonzero and has absolute value less than $p$;
all subsequent terms have larger $p$-adic valuation.
The difference therefore has valuation exactly $\nu_i-1$.
This proves that the sum is direct, has $p^k$ terms, and has the
same branching levels as $A$. It tiles by Newman's criterion.

Two sums from $A^*$ agree modulo $Q_i$ exactly when all their
earlier digits agree. Their difference is then divisible by $g_i$,
and changing digit $i$ alone by one realizes the difference $g_i$.
Thus $\Delta_i(A^*)=\Delta_i$ for every $i$.
Proposition~\ref{prop:digit-classification} gives the same
tilings for $A$ and $A^*$, and also proves that equal data
imply tiling equivalence.

For the converse, let two tiles have the same tiling family.
Choose a common tiling with a period $N$ and residue set $B$.
The identity $|A||B|=N$ gives the same cardinality, hence the same
$p,k$, for the two tiles. Theorem~\ref{thm:periods} recovers $P$
as the smallest fundamental period. A common tiling of period $P$
has branching positions complementary to $I$, by
Lemma~\ref{lem:pure}; hence it recovers all $\nu_i$.
The largest fundamental period recovers $D$.

It remains to recover the lower gcds. Define $w_i=p^{\nu_i-i}$.
For every prime power $q^a\mid D$, formula \eqref{eq:all-counts}
gives
\begin{equation}\label{eq:recover-gcd}
\frac{\log_p F_A(Pq^a)-\log_p F_A(Pq^{a-1})}
{q^a-q^{a-1}}
=\sum_{\substack{1\leq i\leq k\\q^a\mid\Delta_i}}w_i.
\end{equation}
By \eqref{eq:chain}, the indices on the right form a terminal
segment. All $w_i$ are positive, so its sum determines that
segment uniquely. Varying $q^a$ recovers every $\Delta_i$.

Finally, fix $i$ and $m\mid\Delta_i$. Let $\xi,\eta$ be primitive
roots of orders $p^{\nu_i},m$, respectively.
In an occupied class $c\bmod Q_i$, all elements of $A$ have the
same residue modulo $m$. Choose one of them, $a_c$.
The uniform branching recalled in Section~\ref{sec:digits} gives
the contribution
\[
p^{k-i}\eta^{a_c}\xi^c
\sum_{j=0}^{p-1}\xi^{jQ_i}=0
\]
to $A(\xi\eta)$. Therefore $\Phi_{p^{\nu_i}m}(x)\mid A(x)$.
Since $p\nmid\Delta_i$,
\[
\Phi_p(x^{Q_i\Delta_i})
=\prod_{m\mid\Delta_i}\Phi_{p^{\nu_i}m}(x)\mid A(x).
\]
The cyclotomic factors contributed by different $i$ have orders
with distinct $p$-adic valuations.
Their product is $A^*(x)$ and divides $A(x)$.
Comparing degrees proves the inequality in
\eqref{eq:canonical-divisibility}. At equality, the two monic
polynomials are equal, so $A=A^*$ under the normalization
$\min A=0$. Conversely, $A^*$ attains that diameter.
The data are fixed throughout the tiling-equivalence class,
which proves the uniqueness assertion.
\end{proof}

\begin{remark}
The tiling counts alone do not determine the tiling family,
even when the cardinality is prescribed. The canonical binary tiles
with level sets $\{1,2,5,6\}$ and $\{2,3,4,6\}$ and with all
$\Delta_i=1$ both have cardinality $16$ and $P=64$.
For both tiles, $F_A(N)=2^{10}$ when $64\mid N$ and is zero otherwise.
But $\{0,4,8,12\}$ is a tiling of the first in $\mathbb{Z}_{64}$,
whereas for the second the two representations $0+4=4+0$ coincide.
This explains why the preceding proof first recovers the levels
from an actual common tiling before using the counts.
\end{remark}

For the running example $A=\{0,3,10,13\}$, the canonical tile is
\[
A^*=\{0,1\}\oplus\{0,10\}=\{0,1,10,11\}.
\]
It has the same $64$ tilings and diameter $11$ instead of $13$.
In contrast, $\{0,5,10,15\}$ has the same $P=4,D=5$ and the same
branching levels, but its gcd chain is $(5,5)$.
Its CRT fibres may be chosen independently, giving $4^5=1024$
tilings. Thus the lower gcd data cannot in general be replaced
by $D$ alone.

The divisibility $A^*(x)\mid A(x)$ is a polynomial assertion.
It does not assert that $A^*$ is an integer direct-sum factor of $A$;
the quotient polynomial need not have nonnegative coefficients.
For the running example, $A(x)=(1-x+x^2)A^*(x)$.

\begin{corollary}\label{cor:bounds}
The common period and the number of tilings satisfy
\begin{equation}\label{eq:bounds}
PD\leq\frac{p}{p-1}\operatorname{diam}(A),
\qquad
F_A(PD)\leq p^{\operatorname{diam}(A)/(p-1)}.
\end{equation}
Equality in either inequality holds if and only if $A$ is a
translate of an arithmetic progression with $p$ terms.
\end{corollary}

\begin{proof}
The diameter inequality in \eqref{eq:canonical-divisibility} yields
\[
PD+\sum_{i<k}p^{\nu_i}\Delta_i
\leq\frac{p}{p-1}\operatorname{diam}(A).
\]
Hence the first bound is strict for $k>1$.
For the second,
\[
\log_p F_A(PD)
=\sum_{i=1}^k\frac{Q_i\Delta_i}{p^{i-1}}
\leq\sum_{i=1}^kQ_i\Delta_i
\leq\frac{\operatorname{diam}(A)}{p-1},
\]
and the first inequality is strict for $k>1$.
For $k=1$, equality in either bound is precisely the equality case
of Theorem~\ref{thm:canonical}, namely a $p$-term progression.
Conversely, such a progression attains both bounds.
\end{proof}

In particular, all tilings of $A$ share a period at most
$2\operatorname{diam}(A)$, since $p/(p-1)\leq2$.

\section{Affine actions and orbit counts}\label{sec:affine}

The Chinese remainder coordinates give the following affine action.
Fix $M=Pd$, with $d\mid D$, and identify a tiling with its
fibres $C_y\subset\mathbb{Z}_P$ over $\mathbb{Z}_d$.
For a unit $u\bmod M$ and a translation $t\bmod M$, write
$u_P,u_d,t_P,t_d$ for their reductions. The image $uB+t$ has fibres
\begin{equation}\label{eq:affine-fibres}
C'_{u_dy+t_d}=u_P C_y+t_P.
\end{equation}
These maps preserve the tiling family with $A$ fixed:
apply \eqref{eq:dilation} to $B$, since $\gcd(u,|B|)=1$, and then
translate. The cyclic automorphism assertion also appears in
Sands \cite[Theorem~2]{Sands}. A unit multiplier preserves the order
of the translation stabilizer, so fundamental periods are preserved.

To compute fixed tilings, we use successive sections of reduction maps.
For $d\mid D$, define
\begin{equation}\label{eq:section-moduli}
N_0=1,\qquad N_i=p^{\nu_i}d_i,\qquad L_i=N_i/p=Q_i d_i
\quad(1\leq i\leq k).
\end{equation}
Thus $N_{i-1}\mid L_i$ and $N_k=Pd$.
Let $\mathbb{F}_p=\mathbb{Z}_p$ as a field. Start with $B_0=\{0\}\subset\mathbb{Z}_1$,
and, given $B_{i-1}$, set
\begin{equation}\label{eq:section-coordinates}
\begin{split}
H_i&=\{c\in\mathbb{Z}_{L_i}:r_{N_{i-1}}(c)\in B_{i-1}\},\\
\alpha_i&:H_i\longrightarrow\mathbb{F}_p,\\
B_i&=\{c+L_i\alpha_i(c):c\in H_i\}\subset\mathbb{Z}_{N_i}.
\end{split}
\end{equation}
Values of $c$ and $\alpha_i(c)$ in the last line are represented by
$0,\ldots,L_i-1$ and $0,\ldots,p-1$, respectively.
Thus $B_i$ is a section of reduction modulo $L_i$ over $H_i$.

\begin{proposition}\label{prop:section-coordinates}
For $d\mid D$, the construction \eqref{eq:section-coordinates}, with
each $\alpha_i$ arbitrary, gives precisely the members $B_k$ of
$\mathcal{C}_A(Pd)$, without repetition. At every stage,
\[
|H_i|=|B_i|=p^{\nu_i-i}d_i.
\]
Every base $H_k$ arising in this construction admits a section whose
periodic extension has fundamental period $Pd$.
\end{proposition}

\begin{proof}
Under the Chinese remainder isomorphisms, reduction from
$\mathbb{Z}_{N_i}$ to $\mathbb{Z}_{L_i}$ is
\[
\mathbb{Z}_{pQ_i}\times\mathbb{Z}_{d_i}
\longrightarrow\mathbb{Z}_{Q_i}\times\mathbb{Z}_{d_i}.
\]
It forgets only digit $\nu_i$ of the first coordinate. The lift
$c+L_i\alpha_i(c)$ has that digit equal to
\[
\left\lfloor\frac{c}{Q_i}\right\rfloor+d_i\alpha_i(c)
\pmod p.
\]
Since $p\nmid d_i$, arbitrary section values are in bijection with
arbitrary choices of this digit.

Inductively, the points of $H_i$ are labelled by
$\mathcal{D}^{J_i}\times\mathbb{Z}_{d_i}$: the earlier selected
digits are determined by the earlier functions, since $d_h\mid d_i$
for $h<i$, and all intervening digits are free. Thus these are exactly
the independent choices in Proposition~\ref{prop:digit-classification}.
This proves exhaustiveness, uniqueness, and the cardinality formula.

For the last assertion, identify $H_k$ with its image in
$\mathbb{Z}_{Q_k}\times\mathbb{Z}_d$. Choose one point of this base.
Select highest digit $1$ over that point and highest digit $0$ over
all other points. This is an allowed section by the first part.
Every period is a multiple of $P$, by Theorem~\ref{thm:periods},
and hence preserves the first CRT coordinate. It must therefore fix
the unique point of the section whose highest digit is $1$.
A translation of $\mathbb{Z}_{Pd}$ fixing a point is the identity,
so the fundamental period is $Pd$.
\end{proof}

For integers $u,t$ with $\gcd(u,Pd)=1$, take $c$ and $\sigma_i(c)$ as
integer representatives in $\{0,\ldots,L_i-1\}$ and define
\begin{equation}\label{eq:section-carry-notation}
\sigma_i(c)=r_{L_i}(uc+t),\qquad
\kappa_i(c)=\frac{uc+t-\sigma_i(c)}{L_i},\qquad
v=r_p(u)\in\mathbb{F}_p^{\times}.
\end{equation}
If $H_i$ is invariant under $\sigma_i$ and
$\gamma=(c_0,\ldots,c_{\ell-1})$ is a cycle, ordered so that
$\sigma_i(c_j)=c_{j+1}$ with indices modulo $\ell$, define
\begin{equation}\label{eq:cycle-notation}
s_i(\gamma)=\sum_{j=0}^{\ell-1}v^{\ell-1-j}\kappa_i(c_j)
\quad\text{in }\mathbb{F}_p,
\qquad
\omega_i(\gamma)=
\bigl|\{z\in\mathbb{F}_p:(1-v^{\ell})z=s_i(\gamma)\}\bigr|.
\end{equation}

\begin{proposition}\label{prop:section-action}
Let $\gcd(u,Pd)=1$. For the affine map $b\mapsto ub+t$ modulo $Pd$,
the transformed
section has base $\sigma_i(H_i)$ and values
\begin{equation}\label{eq:section-action}
\alpha_i'(\sigma_i(c))=v\alpha_i(c)+\kappa_i(c)
\quad\text{in }\mathbb{F}_p.
\end{equation}
In particular, over a lower set $B_{i-1}$ fixed by the induced affine map,
the number of fixed
sections is the product of $\omega_i(\gamma)$ over the cycles of
$\sigma_i$ on $H_i$, where
\begin{equation}\label{eq:cycle-weight}
\omega_i(\gamma)=
\begin{cases}
1,&v^{\ell}\neq1,\\
p,&v^{\ell}=1\text{ and }s_i(\gamma)=0,\\
0,&v^{\ell}=1\text{ and }s_i(\gamma)\neq0.
\end{cases}
\end{equation}
\end{proposition}

\begin{proof}
For $b=c+L_i\alpha_i(c)$, one has
\[
ub+t=\sigma_i(c)+L_i\bigl(u\alpha_i(c)+\kappa_i(c)\bigr).
\]
Reduction modulo $N_i=pL_i$ proves \eqref{eq:section-action}.
Each $N_i$ divides $Pd$, so reduction commutes with the affine map
at every stage. Changing the representatives of $u,t$ modulo $Pd$
changes $\kappa_i(c)$ by a multiple of $p$, and does not affect the
formula.

If $B_{i-1}$ is fixed, its full inverse image $H_i$ is invariant.
On a cycle $\gamma$, the successive section values are determined
by the initial value $z=\alpha_i(c_0)$. After one circuit the
condition for a fixed section is
\[
(1-v^{\ell})z=s_i(\gamma).
\]
This equation over $\mathbb{F}_p$ has exactly the number of solutions
in \eqref{eq:cycle-weight}. Distinct cycles have independent choices,
which proves the product formula.
The weight $\omega_i(\gamma)$ is independent of the starting point:
changing it transports the initial value by an invertible affine map
of $\mathbb{F}_p$.
\end{proof}

For $N\geq1$ and $\gcd(u,N)=1$, let $F_A(N;u,t)$ and $E_A(N;u,t)$
count the members of $\mathcal{C}_A(N)$ and $\mathcal{E}_A(N)$,
respectively, fixed by $B\mapsto uB+t$.

Starting with $B_0=\{0\}$, retain the fixed sections at each stage
and sum over the fixed lower sets. The last stage gives
$F_A(Pd;u,t)$. If further stages remain, the actual lower sets must
be retained; their numbers alone do not determine the cycle
decomposition at the next stage. This gives a finite recursion for
the fixed-point counts in Proposition~\ref{prop:orbits}.

Write $\mathbb{Z}_N^\times$ for the units modulo $N$. All affine
classes can also be computed at the common
modulus $PD$. For each $d\mid D$, reduction
$\mathbb{Z}_{PD}^{\times}\longrightarrow\mathbb{Z}_{Pd}^{\times}$
is surjective, by the Chinese remainder theorem. Periodic extension
therefore identifies the induced cyclic affine relation on the
stratum of fundamental period $Pd$ with its intrinsic relation modulo $Pd$.
The same assertion holds for translations and signs, and distinct
strata of fixed fundamental period never merge.

Let $R_{\mathrm{tr}}(M)$, $R_{\pm}(M)$, and $R_{\mathrm{aff}}(M)$ denote
the numbers of orbits in $\mathcal{E}_A(M)$ under translations, the
maps $B\mapsto\pm B+t$, and all affine maps $B\mapsto uB+t$ with
$u\in\mathbb{Z}_M^\times$, respectively.

\begin{proposition}\label{prop:orbits}
For $d\mid D$, $M=Pd$ and $\gcd(u,M)=1$,
\begin{equation}\label{eq:fixed-inversion}
E_A(M;u,t)=\sum_{e\mid d}\mu(d/e)F_A(Pe;u,t).
\end{equation}
Moreover,
\begin{align}
R_{\mathrm{tr}}(M)&=\frac{E_A(M)}{M},\label{eq:translation-orbits}\\
R_{\pm}(M)&=\frac{E_A(M)+\sum_{t=0}^{M-1}E_A(M;-1,t)}{2M},
\label{eq:sign-orbits}\\
R_{\mathrm{aff}}(M)&=\frac{1}{M\varphi(M)}
\sum_{u\in\mathbb{Z}_M^\times}\sum_{t=0}^{M-1}E_A(M;u,t).
\label{eq:unit-orbits}
\end{align}
\end{proposition}

\begin{proof}
Periodic extension intertwines the affine map at modulus $M$ with
its reduction modulo $Pe$. Partitioning fixed tilings by
their fundamental period gives
$F_A(M;u,t)=\sum_{e\mid d}E_A(Pe;u,t)$, hence
\eqref{eq:fixed-inversion}.
Translation is free on $\mathcal{E}_A(M)$.
The remaining formulas are Burnside's lemma for the sign/translation
group and the full cyclic affine group. Faithfulness is not needed,
so the sign formula also covers $M=2$.
\end{proof}

Alternatively, the fixed-point counts can be evaluated directly from
the functions in Proposition~\ref{prop:digit-classification} by testing
\eqref{eq:affine-fibres}.

For the running example, the $60$ tilings of fundamental period $20$
give three translation classes, represented by
\[
\{0,2,4,6,8\},\qquad
\{0,4,6,8,12\},\qquad
\{0,2,6,8,14\}.
\]
Each set consists of even residues and meets every residue class
modulo $5$ once, so it is a tiling by \eqref{eq:running-example}.
Their reductions modulo $4$ are not singletons, so each has fundamental period $20$.
Their cyclic gap multisets are, respectively,
$(2,2,2,2,12)$, $(2,2,4,4,8)$, and $(2,2,4,6,6)$,
so these classes are distinct even under reflection.
Writing the three sets as $B^{(1)},B^{(2)},B^{(3)}$, direct calculation gives
\begin{gather*}
-B^{(1)}=B^{(1)}+12,\qquad -B^{(2)}=B^{(2)}+8,\qquad -B^{(3)}=B^{(3)}+12,\\
B^{(3)}=3B^{(1)}+2,\qquad 3B^{(2)}=B^{(2)}+12.
\end{gather*}
Since $3$ and $-1$ generate $\mathbb{Z}_{20}^{\times}$, there are
three sign/translation classes and two cyclic affine classes.
The ordinary affine automorphisms of $\mathbb{Z}$ have multiplier
$1$ or $-1$; general unit multipliers here act on residue sets
before periodic extension.

\subsection*{Declaration of AI use}
The authors were classmates in the Chern class at Nankai University. Starting in December 2024, they held weekly online meetings trying to classify some integer tilings. They are grateful to Terry Tao for his notes~\cite{Tao}, which introduced them to this subject. SageMath, MATLAB, and language models including DeepSeek and Doubao were used to write code for case checks, which helped the authors to observe the structure of fundamental periods, the difference between tile partitions and tile complements, and the normalized tile structure. However the proofs had not been completed until OpenAI's ChatGPT 6 assisted the authors with the project. The authors take full responsibility for the mathematical content of this work. 

\subsection*{Code availability}
\url{https://github.com/LightLearnInTheCampus/matlab_Prime_Power_Tilings_tile_partitions}

\Needspace{24\baselineskip}

\end{document}